\documentclass[twoside,a4paper,reqno,11pt]{amsart}
\usepackage{amsfonts, amsbsy, amsmath, amssymb, latexsym}
\usepackage{mathrsfs,array}
\usepackage[top=30mm,right=30mm,bottom=30mm,left=30mm]{geometry}
\usepackage{stmaryrd}
\usepackage{bm}
\usepackage{xcolor}
\usepackage{hyperref}
\usepackage{tikz-cd}
\usepackage{relsize}

\makeatletter
\@namedef{subjclassname@2020}{\textup{2020} Mathematics Subject Classification}
\makeatother

\def\thefootnote{\fnsymbol{footnote}}

\newcommand{\ovl}{\overline}

\newcommand{\R}{\mathbb{R}}

 \renewcommand{\to}{\rightarrow}

\DeclareMathOperator{\Gal}{Gal}

\DeclareMathOperator{\GL}{GL}

\DeclareMathOperator{\SL}{SL}

\DeclareMathOperator{\red}{red}

\DeclareMathOperator{\Spec}{Spec}

\DeclareMathOperator{\Lie}{Lie}

\DeclareMathOperator{\sep}{s}

\DeclareMathOperator{\Rad}{Rad}
\DeclareMathOperator{\SU}{SU}
\DeclareMathOperator{\Sp}{Sp}
\DeclareMathOperator{\diag}{diag}

\makeatletter
\newcommand{\imod}[1]{\allowbreak\mkern4mu({\operator@font mod}\,\,#1)}
\makeatother

\newtheorem*{conj*}{Conjecture}

\newtheorem{thm}{Theorem}[section]
\newtheorem{prop}[thm]{Proposition}
\newtheorem{lem}[thm]{Lemma}
\newtheorem{cor}[thm]{Corollary}

\theoremstyle{definition}

\newtheorem{example}[thm]{Example}

\begin{document}

\title{Unipotent normal subgroups of algebraic groups}
\author[D.\ Sercombe]{Damian Sercombe}
\address
{Mathematisches Institut, Universität Freiburg, Ernst-Zermelo-Straße 1, 79104 Freiburg, Germany}
\email{damian.sercombe@math.uni-freiburg.de}

\begin{abstract} Let $G$ be an affine algebraic group scheme over a field $k$. We show there exists a unipotent normal subgroup of $G$ which contains all other such subgroups; we call it the restricted unipotent radical $\Rad_u(G)$ of $G$. We investigate some properties of $\Rad_u(G)$, and study those $G$ for which $\Rad_u(G)$ is trivial. In particular, we relate these notions to their well-known analogues for smooth connected affine $k$-groups.
\end{abstract}

\let\thefootnote\relax

\maketitle

\section{Introduction}

\noindent The notion of unipotency plays a key role in the theory of algebraic groups over a field. An algebraic group over a field is \textit{unipotent} if it embeds in the algebraic group of upper unitriangular matrices of $\GL_n$ for some integer $n\geq 1$.

\vspace{2mm}\noindent Over an algebraically closed field, a reductive group is a smooth connected affine group which contains no non-trivial smooth connected unipotent normal subgroups. The theory of reductive groups dates back (at least) to the 1950s; with foundational works by Chevalley and Borel. Reductive groups have a wide range of applications, for example to the Langlands program. 

\vspace{2mm}\noindent Pseudo-reductive groups are an important generalisation of reductive groups, which are relevant only when our underlying field is imperfect. These were first introduced by Borel-Tits in 1978 \cite{BoTi}, and have recently received a lot of attention due to the influential work by Conrad-Gabber-Prasad \cite{CGP} in the 2010s. 
Notably, in 2012 Conrad \cite{Co} used pseudo-reductive groups to prove some finiteness theorems for algebraic groups over global function fields.

\vspace{2mm}\noindent The notions of reductive and pseudo-reductive groups make sense due to the following simple fact. Any affine algebraic group $G$ over a field contains a largest smooth connected unipotent normal subgroup; this is called the \textit{unipotent radical} $\mathscr{R}_u(G)$ of $G$.

\vspace{2mm}\noindent We will show that an analogue of this fact holds without any smoothness or connectedness assumptions. More precisely, we prove the following:

\begin{thm}\label{existproof} Let $k$ be a field. Let $G$ be an affine algebraic $k$-group. There exists a largest unipotent normal subgroup $\Rad_u(G)$ of $G$ (that is, a unipotent normal subgroup of $G$ which contains all other unipotent normal subgroups of $G$). The formation of $\Rad_u(G)$ commutes with base change by separable algebraic field extensions.
\end{thm}

\noindent Henceforth we use the following setup. Let $k$ be an arbitrary field, say of characteristic $p \geq 0$. Let $\overline{k}$ be an algebraic closure of $k$, and let $k^{\sep}$ be the separable closure of $k$ in $\overline{k}$. Let $G$ be an affine algebraic $k$-group (i.e. an affine group scheme of finite type over $k$).

\vspace{2mm}\noindent Let $\Rad_u(G)$ be the largest unipotent normal subgroup of $G$ as in Theorem \ref{existproof}; we call it the \textit{restricted unipotent radical} of $G$. In analogue with the notion of pseudo-reductivity, we say that $G$ is \textit{pseudo-$p$-reductive} if $\Rad_u(G)$ is trivial. Similarly, $G$ is \textit{$p$-reductive} if $\Rad_u(G_{\overline{k}})$ is trivial. We will see that there are parallels between the behaviour of such groups and their smooth connected analogues. However, certain complexities arise when working with $\Rad_u(G)$ that one does not see with $\mathscr{R}_u(G)$. 

\vspace{2mm}\noindent An example of a commutative infinitesimal $k$-group that is pseudo-$p$-reductive but not $p$-reductive is the $k$-group $N$ constructed in Example \ref{tamenoncomm}. 

\vspace{2mm}\noindent Next, we find sufficient conditions on $G$ under which $\Rad_u(G)$ contains \textit{all} unipotent subgroups of $G$ (not just all unipotent normal subgroups). 
Of course this holds if $G$ is commutative. 

\begin{thm}\label{largestunipsub} Let $G$ be an affine algebraic $k$-group. Suppose at least one of the following holds:

\vspace{1mm}\noindent (a) $G$ is trigonalisable over $\overline{k}$;

\vspace{1mm}\noindent (b) $G \cong \overline{G}/Z(\overline{G})$, for some nilpotent connected algebraic $k$-group $\overline{G}$. 

\vspace{2mm}\noindent Then $\Rad_u(G)$ contains all unipotent subgroups of $G$.
\end{thm}


\noindent For example, an affine algebraic $k$-group $G$ is trigonalisable over $\overline{k}$ if it is commutative, 
or if it is smooth, connected and solvable by the Lie-Kolchin theorem \cite[Thm.\ 16.30]{Mi}.

\vspace{2mm}\noindent An example of a nilpotent infinitesimal $k$-group which does not contain a largest unipotent subgroup is the Frobenius kernel $G=(\SL_2)_1$ of $\SL_2$ in characteristic 2. To see this, let $T$ be a maximal torus of $\SL_2$ with corresponding root groups $U^+$ and $U^-$. The center $Z(G)$ of $G$ equals $T_1$, which is of multiplicative type, so $G$ is nilpotent but not unipotent. However $G$ is generated by $(U^+)_1$ and $(U^-)_1$, both of which are 
unipotent. 

\vspace{2mm}\noindent If $G$ is an affine algebraic $k$-group such that all unipotent subgroups of $G$ are contained in $\Rad_u(G)$, then $G/\Rad_u(G)$ contains no non-trivial unipotent subgroups. This motivates the following question: \textit{which affine algebraic $k$-groups contain no non-trivial unipotent subgroups?} 

\vspace{2mm}\noindent Clearly any algebraic $k$-group of multiplicative type satisfies this property; for example tori, or the $n$'th roots of unity $\mu_n$ for any integer $n$. 
A wider class of examples may be constructed as follows: let $k'/k$ be a finite field extension, and let $Z'$ be an algebraic $k'$-group of multiplicative type. Then the Weil restriction $R_{k'/k}(Z')$ is commutative and contains no non-trivial unipotent subgroups.

\vspace{2mm}\noindent It turns out that in many important situations, affine algebraic $k$-groups which contain no non-trivial unipotent subgroups are automatically commutative (in the language of \cite[\S 1.4]{CP}, such groups are called \textit{tame}). More precisely, the following is true. 

\begin{thm}\label{tametheorem} Let $G$ be a connected affine algebraic $k$-group. Suppose $G$ contains no non-trivial unipotent subgroups. Then $G$ is commutative in each of the following cases:

\vspace{1mm}\noindent (a) $k$ is algebraically closed;

\vspace{1mm}\noindent (b) $k$ is perfect and $G$ is solvable; 

\vspace{1mm}\noindent (c) $G$ is trigonalisable over $\overline{k}$.
\end{thm}




\noindent We illustrate Theorems \ref{largestunipsub} and \ref{tametheorem} with several examples of affine algebraic $k$-groups $G$ which contain no non-trivial unipotent subgroups, and yet are neither nilpotent nor trigonalisable over $\overline{k}$. Firstly, over any field $k$, take any perfect finite $k$-group $G$ with order coprime to $p$. 
Secondly, if $k=\R$, take any anisotropic reductive $k$-group $G$ that is not a torus (e.g. $\SU_2$). 
Finally, over an imperfect field $k$, we construct a solvable infinitesimal such $G$ in Example \ref{tamenoncomm}. 

\vspace{2mm}\noindent We now present some consequences of Theorems \ref{largestunipsub} and \ref{tametheorem}. 

\vspace{2mm}\noindent We need the following definition; which is taken from \cite[(1.3.1)]{CGP} (it is most useful for smooth $G$, but makes sense for arbitrary $G$). Let $k'/k$ be the minimal field of definition for $\mathscr{R}_u(\smash{G_{\ovl{k}}})$, and let $\pi':G_{k'} \to G_{k'}/\mathscr{R}_u(G_{k'})=:G'$ be the natural projection. Let \begin{equation}\label{iGmap} i_G:G \to R_{k'/k}(G')\end{equation} be the $k$-homomorphism associated to $\pi'$ under the adjunction of extension of scalars and Weil restriction (refer to \cite[(A.5.1)]{CGP} for this adjunction).

\begin{cor}\label{largestunip} Let $G$ be a solvable smooth connected affine $k$-group. Then $\Rad_u(G)$ contains all unipotent subgroups of $G$. Moreover, $\Rad_u(G)=\ker i_G$, where $i_G$ is the map defined in $(\ref{iGmap})$.
\end{cor}

\begin{cor}\label{largestunipsubcor} Let $G$ be a connected algebraic $k$-group which satisfies either condition (a) or (b) of Theorem \ref{largestunipsub}. Then $G/\Rad_u(G)$ embeds in a Weil restriction $R_{k'/k}(T')$, for some purely inseparable finite field extension $k'/k$ and $k'$-torus $T'$.
\end{cor}

\begin{cor}\label{tametheoremcor} Assume $p>0$. Let $G$ be a connected affine algebraic $k$-group. Suppose the Frobenius kernel $G_1$ of $G$ is of multiplicative type. Then $G$ is of multiplicative type.
\end{cor}

\begin{cor}\label{kakaka} Assume $p>0$. Let $G$ be a connected affine algebraic $k$-group. Suppose the Frobenius kernel $G_1$ of $G$ is trigonalisable over $\overline{k}$ and $p$-reductive. Then $G$ is of multiplicative type.
\end{cor}

\noindent Next, we look at the relationship between pseudo-$p$-reductivity and the usual notion of pseudo-reductivity for smooth connected affine $k$-groups. 

\vspace{2mm}\noindent It is known that, as long as $p>2$, a smooth connected affine $k$-group $G$ is reductive if and only if its Frobenius kernel $G_1$ is $p$-reductive; see for instance \cite[Prop.\ 11.8]{Hu}. The following theorem generalises this result. 

\begin{thm}\label{pseudopreductive} Let $G$ be a smooth connected affine $k$-group. Let $G_1$ be the Frobenius kernel of $G$, let $k'$, $G'$ and $i_G$ be defined as in $(\ref{iGmap})$, and let $k'^{\sep}$ be a separable closure of $k'$.

\vspace{1mm}\noindent (i) Assume $p>2$. Then $G_1$ is pseudo-$p$-reductive if and only if $\ker i_G$ is \'etale. 
If in addition $G$ is perfect, this condition holds if and only if $G$ is pseudo-reductive of minimal type. 

\vspace{1mm}\noindent (ii) Assume $p=2$. Then $G_1$ is pseudo-$p$-reductive if and only if $\ker i_G$ is \'etale and $(G')_{k'^{\sep}}$ has no direct factor isomorphic to $\mathrm{SO}_{2n+1}$ for some $n \geq 1$.
\end{thm}

\begin{cor}\label{pseudopreductivecor} Assume $p>0$. Let $G$ be a smooth connected affine $k$-group. Suppose the Frobenius kernel $G_1$ of $G$ is pseudo-$p$-reductive. Then $G$ is pseudo-reductive and its $k^{\sep}$-root system is reduced.
\end{cor}

\noindent Despite these similarities, the behaviour of $p$-reductive and pseudo-$p$-reductive $k$-groups differs in several important ways from their smooth connected analogues. Notably, $p$-reductivity is not preserved by taking normal subgroups nor by quotients; again refer to Example \ref{tamenoncomm}. 

\section{The restricted unipotent radical}

\noindent In this section we prove Theorem \ref{existproof}. The author thanks Sean Cotner and Timm Peerenboom for some helpful discussions, and the proof of \cite[Prop.\ 5.1.2]{CP} for some ideas.

\vspace{2mm}\noindent Recall that $k$ is an arbitrary field. By a \textit{downward} (resp. \textit{upward}) \textit{directed set}, we mean a partially ordered set in which every finite subset has a lower (resp. upper) bound.

\begin{lem}\label{linalglem} Let $A$ be a $k$-algebra. Let $\{I_i\hspace{0.5mm}|\hspace{0.5mm} i \in \Lambda\}$ be a (possibly uncountable) set of ideals of $A$ that is downward directed under inclusion (i.e. $\Lambda$ is a downward directed set and, for $i,j\in \Lambda$, if $i\leq j$ then $I_i\subseteq I_j$). Let $I:=\bigcap_{i \in \Lambda}I_i$. Then $$A \otimes_k I + I \otimes_k A = \mathsmaller{\bigcap}_{i \in \Lambda}(A \otimes_k I_i + I_i \otimes_k A).$$
\begin{proof} We first check that \begin{equation}\label{firstone} A \otimes_k I + I \otimes_k A = \mathsmaller{\bigcap}_{i \in \Lambda}(A \otimes_k I_i + I \otimes_k A).\end{equation} One inclusion is obvious, for the other let $x \in \bigcap_{i \in \Lambda}(A \otimes_k I_i + I \otimes_k A)$. Take any subspace $V$ of $A$ that complements $I$, and choose a basis for it $\{v_j \hspace{0.5mm}|\hspace{0.5mm} j \in J\}$. Since $x \in (V \otimes_k A) \oplus (I \otimes_k A)$, we can write $x=y+\sum_{j \in J} v_j \otimes a_j$ for unique elements $y \in I \otimes_k A$ and $a_j \in A$. But by assumption $x \in (V \otimes_k I_i) \oplus (I \otimes_k A)$ for every $i \in \Lambda$, hence each $a_j \in \bigcap_{i \in \Lambda}I_i=I$ by uniqueness. So indeed $(\ref{firstone})$ holds.

\vspace{2mm}\noindent Applying $(\ref{firstone})$ twice gives us \begin{equation}\label{secondone}A \otimes_k I + I \otimes_k A = \mathsmaller{\bigcap}_{i,j \in \Lambda}(A \otimes_k I_i + I_j \otimes_k A).\end{equation} We then use the fact that $\Lambda$ is a directed set, giving us \begin{equation}\label{thirdone} \mathsmaller{\bigcap}_{i,j \in \Lambda}(A \otimes_k I_i + I_j \otimes_k A) = \mathsmaller{\bigcap}_{l \in \Lambda}(A \otimes_k I_l + I_l \otimes_k A).\end{equation} 
Combining $(\ref{secondone})$ and $(\ref{thirdone})$ gives us the result. 
\end{proof} 
\end{lem} 

\vspace{2mm}\noindent Now let $X$ be an affine algebraic scheme over $k$. Let $\{Y_i \hspace{0.5mm}|\hspace{0.5mm}i\in \Lambda\}$ be a set of closed subschemes of $X$. 
We define the \textit{schematic union} $\bigcup_{i \in \Lambda}Y_i=:Y$ as follows. Let $\xi:\coprod_{i \in \Lambda} Y_i \to X$ be the morphism of $k$-schemes given by inclusion on each component (indeed it exists, as the category of $k$-schemes admits all -- possibly uncountable -- coproducts). We define $Y$ to be the schematic image of $\xi$, in the sense of \cite[01R7]{Stacks}. Equivalently, $Y$ is the closed subscheme of $X$ whose associated ideal is the intersection of those of all of the $Y_i$'s. 

\vspace{2mm}\noindent In this paper a \textit{subgroup scheme} (or more simply \textit{subgroup}) refers to (the equivalence class associated to) a locally closed immersion of group schemes. Following this convention, by \cite[047T]{Stacks}, all subgroup schemes of a group scheme over a field $k$ are automatically closed.

\begin{prop}\label{schematicim} Let $G$ be an affine algebraic $k$-group. Let $\{H_i\hspace{0.5mm}|\hspace{0.5mm}i \in \Lambda\}$ be a set of (closed) subgroups of $G$ that is upward directed under inclusion. The schematic union $\bigcup_{i \in \Lambda}H_i=:H$ is a (closed) subgroup of $G$. If each $H_i$ is a normal subgroup of $G$ then $H$ is also a normal subgroup of $G$.
\begin{proof} By definition of a $k$-group scheme, $G$ comes equipped with the following morphisms of $k$-schemes: $e:1 \to G$ is the unit map, $m:G \times G \to G$ is multiplication, and $\iota:G \to G$ is inversion. A subscheme $H$ of $G$ is a subgroup if and only if it is stabilised by inversion, and the restriction $m|_{H \times H}$ factors through the inclusion $H \hookrightarrow G$. It is clear that inversion stabilises $H$, since it is an isomorphism of schemes which stabilises each $H_i$, and $H$ is the smallest closed subscheme of $G$ which contains all $H_i$'s. 
We now check multiplication.

\vspace{2mm}\noindent Let $G=\Spec A$ and $H_i=\Spec A/I_i$, for $A$ a finitely generated $k$-algebra and each $I_i$ an ideal of $A$. Denote $I:=\bigcap_{i \in \Lambda}I_i$. By definition $H=\Spec A/I$. Let $m=\Spec \overline{m}$, where $\overline{m}:A \to A \otimes_k A$ is a $k$-algebra homomorphism.

\vspace{2mm}\noindent Fix some $i \in \Lambda$. Since $H_i$ is a subgroup scheme of $G$, the map $$A \xrightarrow{\overline{m}} A \otimes_k A \to A/I_i \otimes_k A/I_i = (A \otimes_k A)/(A \otimes_k I_i + I_i \otimes_k A)$$ 
factors (on the left) through the natural projection $A \to A/I_i$. For $a_i \in I_i$, this means that $$\overline{m}(a_i) \in A \otimes_k I_i + I_i \otimes_k A \subseteq A \otimes_k A.$$ Hence, for $a \in I$, we have that $$\overline{m}(a) \in \mathsmaller{\bigcap}_{i \in \Lambda}(A \otimes_k I_i + I_i \otimes_k A)=A \otimes_k I + I \otimes_k A$$ where the second equality is Lemma \ref{linalglem}. We deduce that the map $$A \xrightarrow{\overline{m}} A \otimes_k A \to A/I \otimes_k A/I = (A \otimes_k A)/(A \otimes_k I + I \otimes_k A)$$ factors (on the left) through the natural projection $A \to A/I$. 

\vspace{2mm}\noindent So indeed $H$ is a subgroup scheme of $G$. It remains to prove the normality assertion.

\vspace{2mm}\noindent We will need the following two facts. Let $\psi:G_1 \to G_2$ be a homomorphism of affine algebraic $k$-groups, and let $\psi^*:\mathcal{O}(G_2) \to \mathcal{O}(G_1)$ be the comorphism (between the coordinate algebras). Firstly, $\psi$ is faithfully flat if and only if $\psi^*$ is injective \cite[Prop.\ 5.43]{Mi}. Secondly, $\mathcal{O}(\ker\psi) = \mathcal{O}(G_1)/\psi^*(I_{G_2})$; where $I_{G_2}=\ker(\mathcal{O}(G_2) \to k)$ is the augmentation ideal \cite[Ch.\ 1e]{Mi}. 

\vspace{2mm}\noindent Now recall our previous setup: $G=\Spec A$, $H_i=\Spec A/I_i$ for each $i \in \Lambda$, and $H=\Spec A/I$ where $I:=\bigcap_{i \in \Lambda}I_i$. Suppose each $H_i$ is a normal subgroup of $G$. For each $i \in \Lambda$, let $f_i:G \to G/H_i=:Q_i$ be the natural projection and say $Q_i= \Spec B_i$. Let $B:=\bigcap_{i \in \Lambda}B_i$ and $Q:=\Spec B$. Choose any $i \in \Lambda$, and define $f:G \to Q$ to be the spectrum of the composite $$B \lhook\joinrel\longrightarrow B_i \xrightarrow{\smash{(f_i)^*}} A.$$ Observe that $f$ is a faithfully flat homomorphism of affine $k$-groups, since $f^*$ is an embedding of Hopf $k$-algebras. 

\vspace{1mm}\noindent Using the aforementioned two facts, we compute \begin{align*}f^*(I_Q) & = f^*(\mathsmaller{\bigcap}_{i \in \Lambda} \ker(B_i \to k))\\ 
& = \mathsmaller{\bigcap}_{i \in \Lambda} (f_i)^*(\ker(B_i \to k)) \\ 
& =  \mathsmaller{\bigcap}_{i \in \Lambda} I_i 
= I \end{align*} Hence $H=\ker f$. This completes the proof. 
\end{proof}
\end{prop}

\noindent We now have enough to prove Theorem \ref{existproof}.

\vspace{2mm}\noindent \underline{Proof of Theorem \ref{existproof}.}

\vspace{2mm}\noindent Let $G$ be an affine algebraic $k$-group. 

\begin{proof} Let $U_1$ and $U_2$ be unipotent normal subgroups of $G$. Since the property of unipotency is preserved by taking quotients and group extensions, $U_1U_2$ is also a unipotent normal subgroup of $G$. So the unipotent normal subgroups of $G$ form an upward directed set under inclusion; let us call it $\mathcal{U}$.

\vspace{2mm}\noindent Let $\Rad_u(G):=\bigcup_{U \in \mathcal{U}} U$; the schematic union. By construction $\Rad_u(G)$ contains all unipotent normal subgroups of $G$. Applying Proposition \ref{schematicim} tells us that $\Rad_u(G)$ is itself a normal subgroup of $G$. We next show that $\Rad_u(G)$ is unipotent. 

\vspace{2mm}\noindent Since $\Rad_u(G)_{\overline{k}} \subseteq \Rad_u(G_{\overline{k}})$, it suffices to assume that $k=\overline{k}$. Consider the (not necessarily Zariski-closed) subgroup $\bigcup_{U \in \mathcal{U}} U(k)$ of the abstract group $G(k)$. 
Since all elements of $\bigcup_{U \in \mathcal{U}} U(k)$ are unipotent, this implies that $\bigcup_{U \in \mathcal{U}} U(k)$ is strictly upper trigonalisable (by the same argument as for closed subgroups). 
Hence $\Rad_u(G)(k)$ is strictly upper trigonalisable, as it equals the Zariski-closure of $\bigcup_{U \in \mathcal{U}} U(k)$. Note that $\Rad_u(G)_{\red}$ is a smooth subgroup of $G$, as $\smash{k=\overline{k}}$. Then applying \cite[XVII, Cor.\ 3.8]{SGA3} tells us that $\Rad_u(G)_{\red}$ is unipotent. This takes care of the case where $p=0$ by Cartier's theorem \cite[Thm.\ 3.23]{Mi}, so we can assume that $p>0$.

\vspace{2mm}\noindent Let $i\geq 1$ be an integer. Consider the schematic union $\bigcup_{U\in\mathcal{U}}U_i$, where $U_i$ denotes the $i$'th Frobenius kernel of $U$. Note that $\bigcup_{U\in\mathcal{U}}U_i$ is a finite $k$-group, as it is contained in $G_i$. 
Hence, since $\mathcal{U}$ is an upward directed set, there exists some
$V\in\mathcal{U}$ such that $V_i=\bigcup_{U\in\mathcal{U}}U_i$. 
Consequently, $\bigcup_{U\in\mathcal{U}}U_i$ is a unipotent normal subgroup of $G$.

\vspace{2mm}\noindent Now choose a sufficiently large integer $t$ such that the chain of subgroups $$\big((\mathsmaller{\bigcup}_{U\in\mathcal{U}}U_i) \cdot \Rad_u(G)_{\red}\big)_{i \geq 1}$$ stabilises for all $i \geq t$ (this is indeed possible as $\Rad_u(G)_{\red}$ is a finite index subgroup of $\Rad_u(G)$). Let $H\in\mathcal{U}$. We apply the following key result taken from \cite[VIIA, Prop.\ 8.3]{SGA3} to $H$: there exists a sufficiently large integer $r$ such that $H/H_r$ is smooth. Since $k=\overline{k}$, this implies that $H=H_r\cdot H_{\red}$. Then $$H
\subseteq (\mathsmaller{\bigcup}_{U\in\mathcal{U}}U_t)\cdot \Rad_u(G)_{\red}.$$ Hence $\Rad_u(G)=(\bigcup_{U\in\mathcal{U}}U_t)\cdot \Rad_u(G)_{\red}$, since $H\in\mathcal{U}$ was chosen arbitrarily. Then, using \cite[Lem.\ 6.41]{Mi}, we deduce that $\Rad_u(G)$ is unipotent.

\vspace{2mm}\noindent The invariance of $\Rad_u(G)$ under base change by separable algebraic field extensions follows from a simple argument of Galois descent. 
Explicitly, let $k'/k$ be a separable algebraic field extension. By uniqueness $\Rad_u(G_{k'})$ is stable under the action of the Galois group $\Gal(k'/k)$, 
so it descends to $k$. It follows that $\Rad_u(G)_{k'} = \Rad_u(G_{k'})$, as the properties of unipotency and normality are invariant under base change. This completes the proof.
\end{proof}

\section{The largest unipotent subgroup}

\noindent Recall that $k$ is a field, of characteristic $p \geq 0$.

\vspace{2mm}\noindent In this section we investigate the conditions under which an affine algebraic $k$-group admits a largest unipotent subgroup. Along the way we prove Theorems \ref{largestunipsub} and \ref{tametheorem}, as well as Corollaries \ref{largestunip}, \ref{largestunipsubcor}, \ref{tametheoremcor} and \ref{kakaka}.

\vspace{2mm}\noindent \underline{Proof of Theorem \ref{largestunipsub}.}

\vspace{2mm}\noindent Let $G$ be an affine algebraic $k$-group (in fact the word affine is superfluous).

\begin{proof} Proof of (a). Suppose $G$ is trigonalisable over $\overline{k}$.

\vspace{2mm}\noindent The formation of the derived group $\mathscr{D}(G)$ commutes with arbitrary base change by \cite[Cor.\ 6.19(a)]{Mi}. Hence $\mathscr{D}(G)$ is unipotent by \cite[Thm.\ 16.6]{Mi}, since $G$ is trigonalisable over $\overline{k}$. Consider the extension $$1 \to \mathscr{D}(G) \to G \xrightarrow{\kappa} G/\mathscr{D}(G) \to 1.$$ Since $G/\mathscr{D}(G)$ is commutative, its restricted unipotent radical $\Rad_u(G/\mathscr{D}(G))$ contains all unipotent subgroups of $G/\mathscr{D}(G)$. In other words, $\Rad_u(G/\mathscr{D}(G))$ is the largest unipotent subgroup of $G/\mathscr{D}(G)$. 
Consider the normal subgroup $U:=\kappa^{-1}(\Rad_u(G/\mathscr{D}(G)))$ of $G$. Recall that the property of unipotency is preserved by taking quotients and group extensions. Since $\mathscr{D}(G)$ is unipotent, it follows that $U$ is the largest unipotent subgroup of $G$. Hence $\Rad_u(G)=U$, by normality of $U$.
\vspace{2mm}\noindent Proof of (b). Let $\overline{G}$ be a connected algebraic $k$-group, and let $G=\overline{G}/Z(\overline{G})$. Then its center $Z(G)$ is unipotent by \cite[XVII, Lem.\ 7.3.2]{SGA3}. [N.B. $G$ is automatically affine by \cite[Cor.\ 8.11]{Mi}.] 

\vspace{2mm}\noindent Suppose in addition that $\overline{G}$, and hence $G$, is nilpotent. Let $n$ be the nilpotency class of $G$. We induct on $n$. If $n=1$ then $G$ is commutative, and the result is given by Theorem \ref{existproof}. 
For $n>1$, consider the central extension \begin{equation}\label{abcde} 1 \to Z(G) \to G \xrightarrow{\zeta} G/Z(G) \to 1.\end{equation}

\vspace{0.5mm}\noindent The nilpotency class of $G/Z(G)$ is $n-1$, so by the inductive hypothesis $\Rad_u(G/Z(G))$ is the largest unipotent subgroup of $G/Z(G)$. Consider the normal subgroup $V:=\zeta^{-1}(\Rad_u(G/Z(G)))$ of $G$. Using the fact that $Z(G)$ is unipotent, a similar argument as in (a) tells us that $\Rad_u(G)=V$.
\end{proof}

\vspace{2mm}\noindent We now move on to the proof of Theorem \ref{tametheorem}. We will repeatedly use the following easy lemma.

\begin{lem}\label{unipcontainsalphap} Assume $p>0$. Let $G$ be an infinitesimal $k$-group. Then either $G$ contains a copy of $\alpha_p$, or $G$ contains no non-trivial unipotent subgroups. 
\begin{proof} Assume that $G$ contains a non-trivial unipotent subgroup $U$. Since $G$ is infinitesimal, then so is $U$. We can then find a copy of $\alpha_p$ in $U$ as the penultimate term of a composition series for $U$; this follows from \cite[XVII, Prop.\ 4.2.1]{SGA3}.
\end{proof}
\end{lem}

\vspace{2mm}\noindent \underline{Proof of Theorem \ref{tametheorem}.}

\vspace{2mm}\noindent Let $G$ be a connected affine algebraic $k$-group which contains no non-trivial unipotent subgroups. 

\begin{proof} Cases (a) and (b). Suppose $k$ is perfect. In this case, showing that $G$ is commutative is equivalent to showing that $G$ is of multiplicative type; this follows from \cite[IV, \S 3, Thm.\ 1.1]{DG}. Suppose in addition that either $k$ is algebraically closed, or $G$ is solvable.

\vspace{2mm}\noindent First, assume that $G$ is smooth. Since $k$ is perfect and $G$ has no non-trivial unipotent subgroups, it is reductive. In the case where $k$ is algebraically closed: $G$ is split and all of its root groups are trivial, so $G$ is a torus. In the case where $G$ is solvable: $G$ is again a torus, as any solvable reductive $k$-group is a torus. So, henceforth, assume that $G$ is not smooth. In particular, $p>0$ by Cartier's theorem \cite[Thm.\ 3.23]{Mi}.

\vspace{2mm}\noindent Next, assume that $G$ is infinitesimal of height $1$. The functor $G \mapsto \Lie(G)$ is an equivalence between the category of height $\leq \!1$ $k$-groups and the category of finite-dimensional restricted Lie algebras; see \cite[II, \S 7, Prop.\ 4.1]{DG}. Under this equivalence, unipotent height $\leq \!1$ $k$-groups correspond to $p$-nilpotent restricted Lie algebras by \cite[XVII, Cor.\ 3.7]{SGA3}. Since $G$ contains no non-trivial unipotent subgroups, $\Lie(G)$ contains no non-trivial $p$-nilpotent elements. We deduce that $\Lie(G)$ -- and hence $G$ -- is commutative (the case where $k$ is algebraically closed is \cite[\S 2, Thm.\ 3.10]{FS}, the case where $G$ is solvable is \cite[\S 2, Prop.\ 3.9]{FS}). Then $G$ of multiplicative type by \cite[IV, \S 3, Thm.\ 1.1]{DG}. 

\vspace{2mm}\noindent Finally, consider the general (non-smooth) case. For any integer $i$, let $G_i$ denote the $i$'th Frobenius kernel of $G$. By \cite[VIIA, Prop.\ 8.3]{SGA3} there exists a sufficiently large integer $r$ such that $G/G_r$ is smooth. Take $r$ to be minimal with respect to this condition, and consider the subnormal series \begin{equation}\label{subnorser} G \supseteq G_r \supset G_{r-1} \supset ... \supset G_1 \supset 1.\end{equation} Each successive quotient of $(\ref{subnorser})$ -- except for the first -- is infinitesimal of height 1. 
By assumption, each subgroup in $(\ref{subnorser})$ contains no non-trivial unipotent subgroups. 

\vspace{2mm}\noindent We induct on the length of this subnormal series $(\ref{subnorser})$. First, assume that the length of $(\ref{subnorser})$ is 1 (that is, $G=G_1$). Then indeed $G$ is of multiplicative type, as we showed previously.

\vspace{2mm}\noindent Next, say $r>1$ and $G=G_r$. Then $G$ is an extension of the height 1 $k$-group $G/G_{r-1}$ by the height $r-1$ $k$-group $G_{r-1}$; let us call the natural projection $\rho_r:G \to G/G_{r-1}$. By the inductive hypothesis, $G_{r-1}$ is of multiplicative type. Assume (for a contradiction) that $G/G_{r-1}$ contains a non-trivial unipotent subgroup. Since $G/G_{r-1}$ is infinitesimal, it contains a subgroup $U$ which is isomorphic to $\alpha_p$ by Lemma \ref{unipcontainsalphap}. Then, since $k$ is perfect and $G_{r-1}$ is of multiplicative type, \cite[XVII, Thm.\ 6.1.1(B)]{SGA3} tells us that the extension $$1 \to G_{r-1} \to \rho_r^{-1}(U) \to U \to 1$$ splits, which is a contradiction. [N.B. In fact $\rho_r^{-1}(U) \cong G_{r-1} \times U$, but we do not need this.] So indeed $G/G_{r-1}$ contains no non-trivial unipotent subgroups. Hence $G/G_{r-1}$ is of multiplicative type by the inductive hypothesis. Then applying \cite[Thm.\ 15.39]{Mi} tells us that $G$ is also of multiplicative type.

\vspace{2mm}\noindent Finally, say $G \supsetneq G_r$. The argument here is very similar to that of the previous paragraph. We have that $G$ is an extension of the smooth $k$-group $G/G_r$ by the infinitesimal $k$-group $G_r$. By the inductive hypothesis, $G_r$ is of multiplicative type. Assume (for a contradiction) that $G/G_r$ contains a non-trivial unipotent subgroup $U$. Then either $U$ is \'etale 
or its Frobenius kernel $U_1$ is non-trivial, in which case $U$ contains a copy of $\alpha_p$ by Lemma \ref{unipcontainsalphap}. In either case \cite[XVII, Thm.\ 6.1.1(A,B)]{SGA3} gives us a contradiction (as \'etale $k$-groups are smooth). Hence $G$ is of multiplicative type, again by \cite[Thm.\ 15.39]{Mi}. In particular, $G$ is commutative. This completes the proof for Cases (a) and (b).

\vspace{2mm}\noindent Case (c). Suppose that $G$ is trigonalisable over $\overline{k}$. Consider the derived group $\mathscr{D}(G)$ as defined in \cite[Def.\ 6.16]{Mi}. Since $G$ is affine, by \cite[Cor.\ 6.19(a)]{Mi} the formation of $\mathscr{D}(G)$ commutes with base change by $\overline{k}/k$. Then base changing by $\overline{k}/k$ and applying \cite[Thm.\ 16.6]{Mi} tells us that $\mathscr{D}(G)$ is unipotent. Hence $\mathscr{D}(G)$ is trivial by assumption. That is, $G$ is commutative. This completes the proof. 
\end{proof}

\noindent \underline{Proof of Corollary \ref{largestunip}.}

\vspace{2mm}\noindent Let $G$ be a solvable smooth connected affine $k$-group. 

\begin{proof} The Lie-Kolchin theorem \cite[Thm.\ 16.30]{Mi} says that $G$ is trigonalisable over $\overline{k}$. Hence $\Rad_u(G)$ contains all unipotent subgroups of $G$, by Theorem \ref{largestunipsub}(a). We now base change by $\overline{k}/k$. A similar argument reassures us that $\Rad_u(G_{\overline{k}})$ contains all unipotent subgroups of $G_{\overline{k}}$. So $G_{\overline{k}}/\Rad_u(G_{\overline{k}})$ contains no non-trivial unipotent subgroups. Then applying Theorem \ref{tametheorem} says that $G_{\overline{k}}/\Rad_u(G_{\overline{k}})$ is commutative, and hence it is of multiplicative type by \cite[IV, \S 3, Thm.\ 1.1]{DG}. Consequently, the extension $$1 \to \Rad_u(G_{\overline{k}}) \to G_{\overline{k}} \to G_{\overline{k}}/\Rad_u(G_{\overline{k}}) \to 1$$ splits by \cite[XVII, Thm.\ 5.1.1]{SGA3}. In particular $\Rad_u(G_{\overline{k}})$ is smooth and connected, so $\Rad_u(G_{\overline{k}})=\mathscr{R}_u(G_{\overline{k}})$. By construction $\ker i_G$ is the largest subgroup of $G$ satisfying the property that $(\ker i_G)_{\overline{k}} \subseteq \mathscr{R}_u(G_{\overline{k}})$. Hence $\ker i_G=\Rad_u(G)$.
\end{proof}

\noindent \underline{Proof of Corollary \ref{largestunipsubcor}.}

\vspace{2mm}\noindent We prove a slightly more general result which implies Corollary \ref{largestunipsubcor}. 

\vspace{2mm}\noindent Let $(P)$ be a property of algebraic group schemes over a field which is preserved by base change by arbitrary (algebraic) field extensions and such that, for any connected affine algebraic $k$-group $G$ with property $(P)$, all unipotent subgroups of $G$ are contained in $\Rad_u(G)$. 

\begin{prop}\label{largestunipsubcorgeneral} Let $G$ be a connected affine algebraic $k$-group $G$ with property $(P)$. Then $G/\Rad_u(G)$ embeds in a Weil restriction $R_{k'/k}(T')$; for some purely inseparable finite field extension $k'/k$ and $k'$-torus $T'$.
\begin{proof} By assumption $G_{\overline{k}}$ also satisfies $(P)$. Consider the restricted unipotent radical $\Rad_u(G_{\overline{k}})$ of $G_{\overline{k}}$. Let $k'/k$ be the minimal field of definition for $\Rad_u(G_{\overline{k}})$; it is purely inseparable as the formation of $\Rad_u(G)$ commutes with separable field extensions, and it is finite by a similar argument as in \cite[Prop.\ 1.1.9(b)]{CGP}. 

\vspace{2mm}\noindent Consider the natural projection $\rho':G_{k'} \to G_{k'}/\Rad_u(G_{k'})=:G'$. Since $\Rad_u(G_{k'})$ is the largest unipotent subgroup of $G_{k'}$ and as the property of unipotency is preserved by group extensions, it follows that $G'$ contains no non-trivial unipotent subgroups. Combining Theorem \ref{tametheorem}(a) with \cite[IV, \S 3, Thm.\ 1.1]{DG} implies that $G'$ is of multiplicative type (as it is a $k'$-descent of a multiplicative type $\overline{k}$-group). Hence $G'$ is contained in some $k'$-torus $T'$. Since the Weil restriction functor preserves embeddings, we have that $R_{k'/k}(G')$ embeds in $R_{k'/k}(T')$.

\vspace{2mm}\noindent Let $\iota:G \to R_{k'/k}(G')$ be the map associated to $\rho'$ under the adjunction of extension of scalars and Weil restriction. The universal property of the counit of this adjunction implies that $(\ker\iota)_{k'} \subseteq \ker\rho'=\Rad_u(G_{k'})$; in particular, $\ker\iota$ is unipotent. So $\ker\iota \subseteq \Rad_u(G)$. In summary, we have embeddings $$G/\Rad_u(G) \hookrightarrow G/\ker\iota \hookrightarrow R_{k'/k}(G') \hookrightarrow R_{k'/k}(T'). \qedhere$$ 
\end{proof} 
\end{prop}

\vspace{1mm}\noindent Clearly properties (a) and (b) of Theorem \ref{largestunipsub} are preserved by base change by arbitrary field extensions. Then Corollary \ref{largestunipsubcor} follows immediately from combining Theorem \ref{largestunipsub} with Proposition \ref{largestunipsubcorgeneral}.

\vspace{4mm}\noindent \underline{Proof of Corollary \ref{tametheoremcor}.}

\vspace{2mm}\noindent Assume $p>0$. Let $G$ be a connected affine algebraic $k$-group. Suppose the Frobenius kernel $G_1$ of $G$ is of multiplicative type.

\begin{proof} The formation of the Frobenius map commutes with arbitrary base change, so without loss of generality we can assume $k=\overline{k}$. If $G$ contains a subgroup $U$ that is isomorphic to $\alpha_p$ then $U$ is contained in $G_1$, which contradicts the fact that $G_1$ is of multiplicative type. So $G$ does not contain a copy of $\alpha_p$. 

\vspace{2mm}\noindent First assume that $G$ is smooth. Since $k$ is algebraically closed and $G$ does not contain a copy of $\alpha_p$, we see that $G$ is reductive and all of its root groups are trivial. Consequently, $G$ is a torus; which is of multiplicative type. 

\vspace{2mm}\noindent Next, assume that $G$ is not smooth. By \cite[VIIA, Prop.\ 8.3]{SGA3} there exists a sufficiently large integer $r$ such that $G/G_r$ is smooth; where $G_r$ denotes the $r$'th Frobenius kernel of $G$. We have a short exact sequence \begin{equation}\label{frobsmooth} 1 \to G_r \to G \xrightarrow{f} G/G_r \to 1.\end{equation} Consider the first Frobenius kernel $(G/G_r)_1$ of $G/G_r$, and its preimage $f^{-1}((G/G_r)_1)=:H$ in $G$. 

\vspace{2mm}\noindent Since $H$ is infinitesimal and does not contain a copy of $\alpha_p$, Lemma \ref{unipcontainsalphap} says that $H$ contains no non-trivial unipotent subgroups. But $k$ is algebraically closed, so applying Case (a) of Theorem \ref{tametheorem} tells us that $H$ is commutative. Hence $H$ is of multiplicative type by \cite[IV, \S 3, Thm.\ 1.1]{DG}. Consequently, its subgroup $G_r$ and its quotient $(G/G_r)_1$ are both also of multiplicative type. We already proved the smooth case, so we deduce that $G/G_r$ is of multiplicative type. Hence $G$ is of multiplicative type, by applying \cite[Thm.\ 15.39]{Mi} to $(\ref{frobsmooth})$.
\end{proof}

\vspace{2mm}\noindent \underline{Proof of Corollary \ref{kakaka}.}

\vspace{2mm}\noindent Assume $p>0$. Let $G$ be a connected affine algebraic $k$-group. Suppose the Frobenius kernel $G_1$ of $G$ is $\overline{k}$-trigonalisable and $p$-reductive.

\begin{proof} Since $G_1$ is $\overline{k}$-trigonalisable and $\Rad_u(G_1)$ is trivial, Theorem \ref{largestunipsub}(a) says that $G_1$ contains no non-trivial unipotent subgroups. Hence $G_1$ is commutative by Theorem \ref{tametheorem}(c). Extending scalars by $\overline{k}/k$ and applying \cite[IV, \S 3, Thm.\ 1.1]{DG} tells us that $G_1$ is of multiplicative type, as it is $p$-reductive. Hence $G$ is of multiplicative type by Corollary \ref{tametheoremcor}. \end{proof}

\vspace{1mm}\noindent In the following example we construct a solvable infinitesimal $k$-group $G$ which is not nilpotent, and contains no non-trivial unipotent subgroups. It is a consequence of Theorem \ref{tametheorem} that $G$ is not trigonalisable over $\overline{k}$.

\begin{example}\label{tamenoncomm} Let $k$ be an imperfect field of characteristic $p>0$, and let $a \in k\setminus k^p$. Let $\mathfrak{g}$ be the 3-dimensional restricted Lie algebra with generators $X,Y,Z$ satisfying $[X,Y]=[X,Z]=0$, $[Z,Y]=Y$, $X^{(p)}=X$, $Y^{(p)}=aX$, $Z^{(p)}=Z$. Explicitly, we have a faithful restricted representation $\mathfrak{g} \to \mathfrak{gl}_p$ given by $X \mapsto \diag(1,1,...,1)$, $Y \mapsto (y_{ij})_{i,j}$ where $y_{1p}=a$, $y_{ij}=1$ if $j=i-1$ and $y_{ij}=0$ otherwise, and $Z \mapsto \diag(1,2,...,p-1,0)$. For instance, if $p=2$ this representation $\mathfrak{g} \to \mathfrak{gl}_2$ is given by $$X \mapsto \begin{pmatrix} 1 & 0 \\ 0 & 1 \end{pmatrix}, ~Y \mapsto \begin{pmatrix} 0 & a \\ 1 & 0 \end{pmatrix}, ~Z \mapsto \begin{pmatrix} 1 & 0 \\ 0 & 0 \end{pmatrix}.$$ 
Our restricted Lie algebra $\mathfrak{g}$ is solvable; it has a unique composition series of restricted subalgebras $$\mathfrak{g} \supset \langle X,Y\rangle \supset X \supset 0.$$ 
Moreover, $\mathfrak{g}$ has no non-zero $p$-nilpotent elements. 

\vspace{2mm}\noindent Let $G:=G(\mathfrak{g})$ be the (unique) height 1 $k$-group associated to the Lie algebra $\mathfrak{g}$ (under the correspondence in \cite[II, \S 7, Prop.\ 4.1]{DG}). Let $N$ be the normal subgroup of $G$ associated to the ideal $\langle X,Y\rangle$ of $\mathfrak{g}$. Then $G$ is a non-nilpotent 
solvable infinitesimal $k$-group of order $p^3$, it has a unique composition series $$G \supset N \supset Z(G) \supset 1$$ with successive quotients $\mu_p$, $\alpha_p$, and $\mu_p$ (taken in this order), and it contains no non-trivial unipotent subgroups. 
Observe that $G$ is $p$-reductive; 
and yet neither its normal subgroup $N$, 
nor its quotient $G/Z(G)$, 
are $p$-reductive. 
\end{example}

\section{Pseudo-$p$-reductivity}

\noindent In this section we prove Theorem \ref{pseudopreductive} and Corollary \ref{pseudopreductivecor}; in which we relate the notion of pseudo-$p$-reductivity with its smooth connected analogue.

\vspace{2mm}\noindent Let $k$ be a field of characteristic $p>0$, and let $G$ be a smooth connected affine $k$-group.

\vspace{2mm}\noindent \underline{Proof of Theorem \ref{pseudopreductive}.}

\vspace{2mm}\noindent Let $k'/k$ be the minimal field of definition for $\mathscr{R}_u(G_{\overline{k}})$, $\pi':G_{k'} \to G_{k'}/\mathscr{R}_u(G_{k'})=:G'$ be the natural projection, and $i_G:G \to R_{k'/k}(G')$ the map associated to $\pi'$ under the adjunction of extension of scalars and Weil restriction.

\begin{proof} We may assume without loss of generality that $k=k^{\sep}$, since the formation of both $\mathscr{R}_u(G)$ and $\Rad_u(G_1)$ commutes with base change by separable algebraic field extensions (by a standard argument of Galois descent). 

\vspace{2mm}\noindent We first prove the ``only if" direction of (i) and (ii).

\vspace{2mm}\noindent Suppose that $G_1$ is pseudo-$p$-reductive. Since $G_1 \cap \ker i_G$ is a unipotent normal subgroup of $G_1$, it must be trivial. Hence $\Lie(\ker i_G)$ is zero, and so $\ker i_G$ is \'etale. Now assume (for a contradiction) that $p=2$ and $G'$ has a direct factor $H'$ that is isomorphic to $\mathrm{SO}_{2n+1}$ for some $n \geq 1$. Consider the unipotent isogeny $$H' \cong \mathrm{SO}_{2n+1} \to \Sp_{2n}$$ described in \cite[\S 2]{PY} or \cite[\S 2]{Va}; its kernel $N'$ is isomorphic to $(\alpha_2)^{2n}$. Since $N'$ is a normal subgroup of a direct factor of $G'$, it is normal in $G'$ itself. 

\vspace{2mm}\noindent Observe that $R_{k'/k}(N')$ is a unipotent normal subgroup of $R_{k'/k}(G')$, since the Weil restriction functor is continuous and it preserves unipotency. Consider the normal subgroup $\smash{V:=i_G^{-1}(R_{k'/k}(N'))}$ of $G$; it is unipotent as the property of unipotency is preserved by taking subgroups and group extensions. The Frobenius kernel $V_1$ of $V$ is a characteristic subgroup of $V$ by \cite[2.28]{Mi}. Consequently, $V_1$ is a unipotent normal subgroup of $G_1$. 
Hence $V_1$ is trivial by assumption, and so $V$ is \'etale. 

\vspace{2mm}\noindent Since $i_G(G)$ is pseudo-split and pseudo-reductive, it contains a Levi subgroup $L$ by \cite[Thm.\ 3.4.6]{CGP}. By \cite[Lem.\ 9.2.1]{CGP} the quotient map $G \to i_G(G)$ induces an isomorphism on the respective maximal reductive $k'$-quotients $\smash{G' \xrightarrow{\sim} i_G(G)'}$. So the inclusion $L \hookrightarrow i_G(G)$ induces an isomorphism $\smash{L_{k'} \xrightarrow{\sim} G'}$. Consequently, since $L$ is split and reductive, it has a direct factor $H$ that is isomorphic to $\mathrm{SO}_{2n+1}$. Let $N \cong (\alpha_2)^{2n}$ be the largest unipotent normal subgroup of $H$. It follows that $R_{k'/k}(N') \cap L$ contains $N$ (one way to justify this is to appeal to \cite[Lem.\ 2.2]{BMRS}). 
But this contradicts the fact that $V$ is \'etale. 

\vspace{2mm}\noindent We next prove the ``if" direction of (i) and (ii).

\vspace{2mm}\noindent Suppose that $G_1$ is not pseudo-$p$-reductive. Recall from \cite[2.27(c)]{Mi} that the formation of the Frobenius kernel commutes with base change by arbitrary field extensions. Since $\pi'$ is smooth and surjective, it restricts to a (faithfully) flat homomorphism $\pi'_1:(G_1)_{k'} \to G'_1$ by \cite[VIIA, Prop.\ 8.2]{SGA3}; where $G'_1$ denotes the Frobenius kernel of $G'$.

\vspace{2mm}\noindent Let $U$ be a non-trivial unipotent normal subgroup of $G_1$. Since $\pi'_1$ is (faithfully) flat, $\pi'(U_{k'})=:U'$ is a unipotent normal subgroup of $G'_1$. There are two possibilities: either $U_{k'} \subseteq \ker\pi'$ or $U'$ is non-trivial.

\vspace{2mm}\noindent First consider the case where $U_{k'} \subseteq \ker\pi'=\mathscr{R}_u(G_{k'})$. By construction $\ker i_G$ is the largest subgroup of $G$ satisfying the property that $(\ker i_G)_{k'} \subseteq \mathscr{R}_u(G_{k'})$. So $U$ is contained in the unipotent normal subgroup $G_1 \cap \ker i_G$ of $G_1$. Hence $\ker i_G$ is not \'etale.

\vspace{2mm}\noindent Next consider the case where $U'$ is non-trivial. Denote $\mathfrak{g}':=\Lie(G')=\Lie(G'_1)$, and consider the ideal $\mathfrak{u}':=\Lie(U')$ of $\mathfrak{g}'$. Since $U'$ is unipotent it intersects trivially with any maximal torus $T'$ of $G'$, hence $\mathfrak{u}' \cap \Lie(T')=0$. Then applying \cite[Lem.\ 2.1]{Va} tells us that $p=2$ and $G'$ has a direct factor which is isomorphic to $\mathrm{SO}_{2n+1}$ for some $n \geq 1$. [N.B. We cannot simply deduce this from \cite[Thm.\ 2.2(a)]{Va} or \cite[Lem.\ 2.2]{PY} as the restricted unipotent radical is not necessarily a characteristic subgroup.] 

\vspace{2mm}\noindent It remains to prove the second assertion of (i). Suppose $p>2$. 

\vspace{2mm}\noindent Assume that $G$ is pseudo-reductive of minimal type. Since $p >2$, the root system of $G$ is reduced by \cite[Thm.\ 2.3.10]{CGP}. Hence $\ker i_G$ is central in $G$ by \cite[Prop.\ 6.2.15]{CP1}. By definition of minimal type, this implies $\ker i_G$ is trivial.

\vspace{2mm}\noindent For the converse, assume that $G$ is perfect and $\ker i_G$ is \'etale. Since $G$ is pseudo-split and pseudo-semisimple, it decomposes as a commuting product $G=S_1...S_r$; where each $S_j$ is a pseudo-split (absolutely) pseudo-simple $k$-group. It follows that $\ker i_G = \ker i_{S_1}...\ker i_{S_r}$; 
in particular, $\ker i_{S_j}$ is \'etale for each $j=1,...,r$. We need the following result taken from \cite[Ex.\ 7.1.1]{CP1}: if $p>2$, $S$ is an absolutely pseudo-simple $k$-group and $\ker i_S$ is \'etale, then $\ker i_S$ is trivial. Applying this to our decomposition says that each $\ker i_{S_j}$ is trivial, hence so is $\ker i_G$. It follows that $G$ is pseudo-reductive of minimal type. This completes the proof.
\end{proof}

\vspace{2mm}\noindent \underline{Proof of Corollary \ref{pseudopreductivecor}.}

\begin{proof} Assume $G_1$ is pseudo-$p$-reductive. Applying Theorem \ref{pseudopreductive} tells us that $\ker i_G$ is \'etale. Then $G$ is pseudo-reductive, as $\mathscr{R}_u(G) \subseteq \ker i_G$. Moreover, $\ker i_G$ is central in $G$ by \cite[Rem.\ 12.39(a)]{Mi}. The central quotient map $G \to i_G(G)$ induces a bijection on the respective root systems by \cite[Prop.\ 2.3.15]{CGP}. But $i_G(G)$ is pseudo-reductive with reduced $k^{\sep}$-root system, so the root system of $G_{k^{\sep}}$ is also reduced.
\end{proof}

\bigskip\noindent
{\textbf {Acknowledgements}}:
The author would like to thank Michael Bate, Brian Conrad, Sean Cotner, Ben Martin, Timm Peerenboom, Gopal Prasad, Gerhard Röhrle and David Stewart, as well as an anonymous referee, for some helpful comments and discussions. Most of this research was funded by an Alexander von Humboldt postdoctoral grant, and currently the author is employed/funded by the University of Freiburg Mathematical Institute.

\bibliographystyle{amsalpha}

\newcommand{\etalchar}[1]{$^{#1}$}
\providecommand{\bysame}{\leavevmode\hbox to3em{\hrulefill}\thinspace}
\providecommand{\MR}{\relax\ifhmode\unskip\space\fi MR }
\providecommand{\MRhref}[2]{%
	\href{http://www.ams.org/mathscinet-getitem?mr=#1}{#2} }
\providecommand{\href}[2]{#2}

\end{document}